\documentclass[11pt]{amsart}
\usepackage{amsfonts}
\usepackage{amsfonts,latexsym,rawfonts,amsmath,amssymb,amsthm}
\usepackage[plainpages=false]{hyperref}

\usepackage{graphicx}

\numberwithin{equation}{section}

\newcommand{\beq}{\begin{equation}}
\newcommand{\eeq}{\end{equation}}
\newcommand{\beqs}{\begin{eqnarray*}}
\newcommand{\eeqs}{\end{eqnarray*}}
\newcommand{\beqn}{\begin{eqnarray}}
\newcommand{\eeqn}{\end{eqnarray}}
\newcommand{\beqa}{\begin{array}}
\newcommand{\eeqa}{\end{array}}

\def\lra{\longrightarrow}

\def\bc{\begin{center}}
\def\ec{\end{center}}

\def\begeq{\begin{equation}}
\def\endeq{\end{equation}}
\def\and{\quad{\rm and}\quad}

\let\lra=\longrightarrow

\def\mapright\#1{\,\smash{\mathop{\lra}\limits^{\#1}}\,}

\newtheorem{prop}{Proposition}[section]
\newtheorem{theo}[prop]{Theorem}
\newtheorem{lem}[prop]{Lemma}

\newtheorem{cor}[prop]{Corollary}
\newtheorem{rem}[prop]{Remark}

\newtheorem{defi}[prop]{Definition}

\title  {Minimizing  weak solutions for calabi's extremal metrics on  toric manifolds}

\author {Bin $\text{Zhou}^*$}
\address{ * Department of Mathematics, Peking University,
Beijing, 100871, China}
\author { Xiaohua $\text{Zhu}^{*,**}$}
\thanks { **  Partially supported by  NSF10425102 in China and the Huo Y-T Fund.}
 \email{xhzhu\@math.pku.edu.cn}

\begin{document}
\bibliographystyle{plain}

\begin{abstract}In this paper, we discuss a Donaldson's version of the modified $K$-energy
associated to the Calabi's extremal metrics on  toric manifolds and
prove the existence of the weak solution for extremal metrics in the
sense of convex functions which minimizes the modified $K$-energy.
\end{abstract}
\maketitle

\setcounter{section}{-1}

\section { Introduction }

The existence of extremal metrics has been recently studied
extensively on K\"ahler manifolds. The goal is to establish a
sufficient and necessary  condition for the existence of extremal
metrics in the sense of Geometric Invariant Theory.  There are many
important works related to the necessary part ([Ti], [D1],
[M1],[M2]). The sufficient part seems more difficult than the
necessary part since the existence is related to the solvability of
certain fourth-order elliptic equations. On the other hand, an
extremal metric can be regarded as a critical point of some
geometric energies, such as the Calabi's energy, the modified
$K$-energy, e.t. This gives a way to study the existence by using
variational method in the sense of Nonlinear Analysis.  In this
paper, we focus on a class of special K\"ahler manifolds, namely,
toric manifolds and discuss the minimizing weak solution for
extremal metrics in the sense of convex functions related to a
Donaldson's version of the modified $K$-energy.

Let $(M, g)$ be a compact K\"ahler manifold of dimension $n$. Then
K\"ahler form $\omega_g$ of $g$ is given  by
$$\omega_g=\sqrt{-1}\sum_{i,j=1}^n g_{i\overline j} dz^i\wedge d\overline
z^j$$
 in local coordinates $(z_1,...,z_n)$, where $g_{i\overline j}$ are components of metric $g$.
 The K\"ahler class $[\omega_g]$ of $\omega_g$ can be represented by a set of potential
 functions as follow
 $$\mathcal {M} =\{\phi\in C^{\infty}(M)|
~\omega_{\phi}=\omega_g+ \frac {\sqrt {-1}}{2\pi}
\partial\overline\partial\phi >0\}.$$
According to [Ca], a K\"ahler metric $\omega_{\phi}$ in the
K\"ahler class $[\omega_g]$ is called extremal if
  \beq R(\omega_{\phi})=\overline R+\theta_X(\omega_\phi)\eeq
   for some holomorphic vector field $X$ on $M$, where
   $R(\omega_{\phi})$ is the
scalar curvature of $\omega_{\phi}$,  $\overline R=\frac
{1}{V}\int_M R(\omega_g)\omega_g^n$, $V=\int_M \omega_g^n$ and
$\theta_X(\omega_\phi)$ denotes the potential function of $X$
associated to the metric $\omega_\phi$,  which is defined by
$$i_X\omega_{\phi} =  \sqrt{-1}\overline\partial \theta_X(\omega_{\phi}),
\int_M \theta_X(\omega_{\phi})\omega_\phi^n=0.
$$
By [FM], such an  $X$, usually called extremal is uniquely
determined by the Futaki invariant $F(\cdot)$.

 Equation (0.1) can be regarded as an Euler-Langrange equation  of
the following modified $K$-energy
$$\mu (\phi) = - \frac{1}{V}\int_{0}^{1}\int_{M}
\dot\psi_{t}[R(\omega_{\psi_{t}}) -\overline R- \theta_X(\psi_t)]
\omega_{\psi_{t}}^n\wedge dt,$$
  where $\psi_{t} (0 \leq t \leq 1)$
is a path connecting $0$ to $\phi$ in $\mathcal {M}$.  In fact, one
can show that the functional $\mu(\phi)$ is well-defined, i.e., it
is independent of the choice of path $\psi_t$ ([Gua]). Thus $\phi$
is a critical point of $\mu(\cdot)$ iff the corresponding metric
$\omega_\phi$ is extremal.

\begin {defi} Let
$$I(\phi)=\frac{1}{V}\int_M \phi(\omega_g^n-\omega_{\phi}^n).$$
$\mu(\phi)$ is called proper associated to a subgroup $G$ of the
automorphisms group $\text{Aut}(M)$ in K\"ahler class $[\omega_g]$
if there is a continuous function $p(t)$ in $\Bbb R$ with the
property

$$\lim_{t\to+\infty} p(t)=+\infty,$$
  such that
  \beq\mu(\phi)\ge \inf_{\sigma\in G} p(I(\phi_{\sigma})),\eeq
   where
 $\phi_{\sigma}$ is defined by
$$\omega_g+\sqrt{-1}\partial\bar{\partial}\phi_{\sigma}
=\sigma^*(\omega_g+\sqrt{-1}\partial\bar{\partial}\phi).$$
\end{defi}

There is a natural question: does there exists an extremal metric if
the modified K-energy $\mu(\phi)$ is proper associated to a
reductive subgroup $G$ of $Aut(M)$? The answer for K\"ahler-Einstein
metrics with positive scalar curvature is positive by Tian  ([Ti]).
The converse is also true. In this paper,  we discuss this question
for toric manifolds.

An n-dimensional toric manifold $M$ corresponds to a polytope $P$ in
$\mathbb {R}^n$ which satisfies Delzant's condition ([Gui]).  Let
$G_0$ be a maximal compact subgroup of  torus actions group $T$ on
$M$. It was showed in [ZZ] that for any $G_0$-invariant $\phi$,
  \beq \mu(\phi)=\frac{2^nn!(2\pi)^n}{V} \mathcal
{F}(u)\eeq
  with
\beq  \mathcal {F}(u)=-\int_P \log(\det(D^2u))dx +\int_{\partial P}u
d\sigma-\int_P (\bar{R}+\theta_X)udx,\eeq
 where $u$ is a Legendre function of $\phi$ which is smooth convex
 function in $P$ and can be extended to a continuous function on $\overline P$,
  and $d\sigma$ is a natural induced measure from $dx$. In [D2],
  Donaldson first obtained (0.3) for the $K$-energy. Since Hessian matrix $(D^2u)$ exists almost
everywhere for a convex function $u$ in $P$, it is possible to
 extend $\mathcal {F}(u)$ to a more  general class of convex functions in $P$.

Let  $\tilde{\mathcal {C}}$ be a set of normalized Legendre
functions associated to $G_0$-invariant potential functions on $M$
(cf. Section 1). Set
$$\mathcal {C}_{\star}=\{u\geq 0 \text{~is a limit of some
sequence of \{$u_n$\} in $\tilde{\mathcal {C}}$ with $\int_{\partial
P}u_nd\sigma<C$}\}.$$
 Then $\mathcal {C}_{\star}$ is a complete space in the sense of local $C^0$-convergence.
Moreover we show that $\mathcal {F}(u)$ is well-defined in $\mathcal
{C}_{\star}$ although it may be infinity (cf. Section 2).  It is
easy to see that the Euler-Langrange equation for $\mathcal {F}(u)$
in $\tilde{\mathcal {C}}$   is
 \beq
-u_{ij}^{ij}=\bar{R}+\theta_X, ~\text{in} ~P\eeq
 where $\theta_X$ is a potential function of  $X$  which is a linear
 function in $P$.  We call $u$ a weak solution
  of (0.5) in the sense of convex functions if  $\mathcal {F}(u)<\infty$ and
  $u$ is a critical point of $ \mathcal {F}(u)$
  in  $\mathcal {C}_{\star}$.

The following is our main theorem in this paper.

\begin {theo} Suppose that $\mu(\phi)$ is proper  for $G_0$-invariant potential functions  associated to toric
actions group  $T$ on a toric manifold $M$. Then there exists a weak
solution $u_\infty$ of equation (0.5) for extremal metrics on $M$ in
the sense of convex functions which minimizes $\mathcal {F}(u)$ in
$\mathcal{C}_{\star}$.
\end {theo}

In [ZZ], authors  introduced a sufficient condition to verify the
properness of $\mu(\phi)$ on toric manifolds and found some examples
which satisfy the condition. Thus according to the above theorem
there exists a minimizing weak solution of (0.5) on these toric
manifolds. The regularity of minimizing solution of (0.5) is an
interested topic. If the minimizing solution $u_\infty$ in Theorem 2
belongs to $\tilde{\mathcal {C}}$, then $u_\infty$ induces an
extremal metric on $M$.  We hope to discuss this problem in the
future.

The organization of this paper is as follows. In Section 1, we
review the modified $K$-energy $\mu(\phi)$ in the sense of
Donaldson's version  for convex functions on toric manifolds. In
section we will  extend   $ \mathcal {F}(u)$ to more  general convex
functions.  In section 3, we prove the lower semi-continuity of $
\mathcal {F}(u)$.   In section 4, we prove  Theorem 0.2 and discuss
some properties about the   minimizing weak solution of (0.5).

\vskip3mm
  \noindent {\bf Acknowledgements} The authors would like to thank
  professor X-J Wang for many valuable discussions.

\section { Modified K-energy  $\mu(\phi)$}
 \vskip3mm

In this section, we recall a Donaldson's version of the modified
K-energy $\mu(\phi)$ of  on toric manifolds. We assume that $M$ is
an n-dimensional toric K\"ahler manifold and $g$  is a $G_0\cong
(S^1)^n$-invariant K\"ahler metric in the K\"ahler class,  where
$G_0$ is a maximal compact subgroup of torus actions group $T$ on
$M$. Then under an affine logarithm coordinates system
$(w_1,...,w_n)$, its K\"ahler form $\omega_g$ is determined by a
convex function $\psi_0$ on $\mathbb {R}^n$, namely,
$$\omega_g = \sqrt{-1}\partial\bar{\partial}\psi_0$$
is defined on the open dense orbit $T$. Denote $D\psi_0$ to be a
 gradient map (moment map) associated to $T$.
 Then the image of $D\psi_0$ is a convex polytope
$P$ in $\mathbb {R}^n$.
 By using the Legendre transformation $y=(D\psi_0)^{-1}(x)$, we see
that the function (Legendre function) defined by
$$u_0(x)=\langle y,D\psi_0(y) \rangle -\psi_0(y)=\langle y(x),x\rangle-\psi_0(y(x))$$
is  convex on $P$.  In general, for any $G_0$-invariant potential
function $\phi$  associated to the K\"ahler class $[\omega_g]$,
one  gets a convex function $u(x)$ on $P$ by using the above
relation  while $\psi_0$ is replaced by $\psi_0+\phi$. Set
  $$\mathcal {C} =\{u=u_0+v| ~ u~ \text{is a convex function in }~ P,~v\in
C^{\infty}(\bar{P})\}.$$
 It was showed in [Ab] that functions in $\mathcal {C}$ are
 corresponding to
$G_0$-invariant functions in $\mathcal {M}$ (whose set is  denoted
by $\mathcal {M}_{G_0}$)  by one-to-one.

An n-dimensional toric manifold $M$ corresponds to a polytope $P$ in
$\mathbb {R}^n$ which is described by a common set of some
half-spaces,
  $$\langle l_i, x\rangle < \lambda_i, ~i=1,...,d, $$
 where $l_i$ are $d-$vectors in $\mathbb {R}^n$ with all components in $\mathbb {Z}$,
 which satisfy the Delzant condition ([Gui]). Without  loss of
 generality, we may assume that the original point $0$ lies in $P$, so all $\lambda_i>0$.
A special element of $\mathcal {C}$ in the sense of $P$ can be
constructed  as follow ([Gui]),
 \beq u_P=\sum (-\langle l_i, x\rangle +
\lambda_i)\log (-\langle l_i, x\rangle + \lambda_i).\eeq
 The convex function $u_P$ is very useful and we will use it at many places in our paper.

In [D2], Donaldson found a formula for the $K$-energy  in the sense
of convex functions in $\mathcal {C}$. For the modified $K$-energy
$\mu(\phi)$, we have ([ZZ]),

\begin {lem} Let $d\sigma_0$ be the Lebesgue measure on the
boundary $\partial P$ and $\nu$ be the outer normal vector field on
$\partial P$. Let $d\sigma=\lambda_i^{-1}(\nu, x)d\sigma_0$ on the
face $\langle l_i, x\rangle= \lambda_i$ of $P$. Then
 \beq
\mu(\phi)=\frac{2^nn!(2\pi)^n}{V} \mathcal {F}(u),\eeq
 where
\beq \mathcal {F}(u)=-\int_P \log(\det(D^2u))dx +\int_{\partial
P}u d\sigma-\int_P (\bar{R}+\theta_X)udx. \eeq
\end {lem}

The functional $\mathcal {F}(u)$ is invariant according to the
choice of $X$ if $u$ is replaced by adding an affine linear
function. For this reason, we normalize $u$ as follows. Let $p\in P$
and set
 $$\tilde{\mathcal {C}}= \{u \in \mathcal {C}
|\inf_{P}u=u(p)=0 \}.$$
 Then for any $u_{\phi} \in \mathcal {C}$ corresponding to $\phi\in \mathcal {M}_{G_0}$,
one can normalize $u_{\phi}$ by
$$\tilde{u}_{\phi} = u_{\phi} - (\langle Du_{\phi}(p),x-p\rangle + u_{\phi}(p))$$
so that $\tilde{u}_{\phi}=u_{\tilde\phi}\in \tilde{\mathcal {C}}$
corresponds to a K\"ahler potential function $\tilde{\phi}\in
\mathcal {M}_{G_0}$  which satisfies
$$D(\tilde{\phi} + \psi_0)(0) = p~\text{and}~ (\tilde\phi+\psi_0)(0)=0.$$
In fact, $\tilde\phi$ can be uniquely determined by using the
affine coordinates transformation $y\to y+y_0$ as follow,
$$\tilde\phi(y)=(\phi+\psi_0)(y+y_0)-\psi_0(y)- (\phi+\psi_0)(y_0).$$

The following lemma was also proved in [ZZ].

\begin {lem} There is a constant $C$ independent of $\phi$, such
that \beq |\int_P\tilde{u}_{\phi}dx-I(\tilde\phi)|\leq C. \eeq
\end {lem}

By Lemma 1.2, one can easily get

\begin {prop} $\mu(\phi)$ is proper for $G_0$-invariant potential functions   associated to toric actions $T$
 if only if there exists  a continuous
function $p(t)$ in $\Bbb R$ with the property
 $$\lim_{t\to+\infty} p(t)=+\infty,$$
  such that
  $$ \mathcal {F}(u)\ge p(\int_P u dx), ~\forall u\in \tilde{\mathcal {C}}.$$
  \end {prop}

At the end of this section, we  would like to recall a set
$$\mathcal {C}_{\infty}=\{u\in C^{\infty}(P)\cap C(\overline
P)|~u~\text{is convex in}~ P\}.$$
 Clearly, $\mathcal {C} \subset \mathcal {C}_{\infty}$. It was proved in [D2]
    $$ \mathcal {F}(u)>-\infty,~\forall ~u\in \mathcal
 {C}_{\infty}$$
   and
 \beq\inf_{\mathcal {C}_{\infty}}\mathcal
 {F}(u) =\inf_{\mathcal {C}}\mathcal
 {F}(u).\eeq

\section {  Lower bound of $\mathcal {F}(u)$}
 \vskip3mm

Since $\tilde{\mathcal {C}}$ is not a complete space, we need to
define a closure set of $\tilde{\mathcal {C}}$.  We denote by
$P^{\ast}$ the union of $P$ and the open codimension-1 faces. And
define
$$\mathcal {C}_{\star}=\{u\geq 0 \text{~is a limit of some sequence of \{$u_n$\} in $\tilde{\mathcal {C}}$ with
$\int_{\partial P}u_nd\sigma_0<C$}\}.$$
  Note that for any
$u\in{\mathcal {C}_{\star}}$, by defining $u$ on the boundary to
be
$$\lim_{t\rightarrow 1}u_{\infty}(tz),$$
then it is a convex function
on $P^{\ast}$ with $\int _{\partial P}ud\sigma_0 < \infty$. On the
other hand, $\int_{\partial P}u_nd\sigma_0<C$ implies that there is
a subsequence of $u_n$ converging uniformly on any compact set of
codimension$-1$ faces. We denote the limit to be $\hat{u}$ on any
open set of  codimension$-1$ faces.  By the convexity of $u$, it is
easy to see that
$$u|_{\partial P}\leq\hat{u}.$$

Let $u_0\in{\tilde{\mathcal {C}}}$ and choose $K_0>0$ so that
$$\int_{\partial P}u_0 d\sigma_0 <K_0.$$
For any $K>K_0$, we denote a subset of convex functions in
$P^{\ast}$ by
 \beqs  \mathcal {C}_{\star}^K(P)&&=\{u \text{~is a
convex function in }~P^{\ast}~ \text{with}~\\
&& \int_{\partial P}u d\sigma_0\leq K~ \text{and}~
\inf_{P}u=u(0)=0\}. \eeqs
  Clearly $u_0\in{\mathcal
{C}_{\star}^K(P)}$. Note that for any sequence $u_n$ converging
locally uniformly to $u$ in $\mathcal {C}_{\star}$,
 $$\int_{\partial P} u d\sigma_0=\int_{\partial P}\hat{u} d\sigma_0=\lim_n\int_{\partial
P}u_nd\sigma_0\leq K.$$
 We see  that all $\mathcal {C}_{\star}^K(P)$
are complete spaces in sense of local $C^0-$convergence.

We want to extend $\mathcal {F}(u)$ to be defined on $\mathcal
{C}_{\star}^K(P)$. Since the linear part
$$L(u)=\int_{\partial P}u
d\sigma-\int_P (\bar{R}+\theta_X)udx$$ is well-defined on $\mathcal
{C}_{\star}^K(P)$, it  suffices to define the nonlinear part
$$\int_P\log(\det(D^2u))dx.$$

Note that $u$ is almost everywhere twice-differentiable since $u$ is
convex and so the Hessian matrix $(D^2u)$ exists almost everywhere.
We denote Hessian matrix $(\partial^2u)$ at those
twice-differentiable points in $P$. Then one can show that
$\det(\partial^2u)$ corresponds to the regular part $\mu_r[u]$ of
Monge-Amp$\grave{e}$re measure $\mu[u]$ associated to $u$ and so it
is a locally integrable function ([TW2]). Thus
 $$\int_P \log^+(\det(\partial^2u))dx$$
  is well-defined although the integral
may be infinity, where
$$\log^+(\det{(\partial^2 u)})=\max(0,
\log(\det(\partial^2 u))).$$
  But we have

\begin {prop}Let $u\in{\mathcal {C}_{\star}^K(P)}$. Then
$$\int_P\log^+(\det(\partial^2u))dx<\infty.$$
\end {prop}

By Proposition 2.1, one sees that
 $$\int_P \log(\det(\partial^2u))dx$$
  is integrable.  To prove the proposition,   we need to regularize $u$ by
mollification functions, namely,
$$u_h(x)=h^{-n}\int_P\rho(\frac{x-y}{h})u(y)dy$$
for any small $h>0$, where $x$  satisfies $h\leq \text{dist}(x,
\partial P)$ and  $\rho$ is a support function in $B_1(0)\subset
\mathbb {R}^n$ with $\int_{B_1(0)}\rho=1$.  A fundamental result is
that $(D^2u_h)\rightarrow (\partial^2u)$ almost everywhere for a
convex function $u$ in $P$.

\begin {lem}Let $u\in {\mathcal {C}_{\star}^K(P)}$ and $u_n$ be
a sequence in $C^2(P)$ converging locally uniformly to $u$ with
$\det(\partial^2u_n)\rightarrow\det(\partial^2u)$ almost
everywhere. Suppose that

\beq \det(\partial^2 u_n),~\det(\partial^2 u)\geq\epsilon_0>0.
\eeq Then for any subset $\Omega\subset\subset P$,
$$\int_{\Omega}\log(\det(\partial^2u))dx=\lim_{n\rightarrow \infty}\int_{\Omega}\log(\det(D^2u_n))dx.$$
\end {lem}

\begin {proof}  Let $\Omega\subset P$ be any Borel subset  with
$\text{dist}(\Omega,\partial P)\geq \delta>0$.  By using the
concavity of $\log$, we have
 \beqs \int_\Omega
\log(\det(\partial^2u))dx
&&\leq |\Omega|\log\left(|\Omega|^{-1}\int_\Omega \det(\partial^2u)dx\right)\\
&&\leq |\Omega|\log\left(|\Omega|^{-1}(\frac{\text{osc}(u)}{\delta})^n\right)\\
&&=n|\Omega|\frac{\text{osc}(u)}{\delta}-|\Omega|\log(|\Omega|).\eeqs
Similarly,  for $u_n$, it holds
 \beq \int_\Omega \log(\det(\partial^2u_n))dx
\leq
n|\Omega|\frac{\text{osc}(u_n)}{\delta}-|\Omega|\log(|\Omega|).\eeq

For any $\lambda>0$,  we let
$$\Omega_{n,\lambda}=\{x\in {\Omega}|~ |\log(\det(\partial^2u_n))-\log(\det(\partial^2u))|>\lambda\}.$$
Then $|\Omega_{n,\lambda}|\rightarrow 0$ as $n\rightarrow \infty$
since $D^2u_n\rightarrow D^2u$ almost everywhere. Note
 \beqs
&&|\int_{\Omega}\log(\det(\partial^2u_n))dx-\int_{\Omega}\log(\det(\partial^2u))dx|\\
&& \leq
|\int_{\Omega_{n,\lambda}}\log(\det(\partial^2u_n))dx|+|\int_{\Omega_{h,\lambda}}\log(\det(\partial^2u))dx|\\
&&+\int_{\Omega\setminus\Omega_{n,\lambda}}|\log(\det(\partial^2u_n))dx-\log(\det(\partial^2u))dx|\eeqs
 It is clear that the first two integrals  in  the right side of  the inequlity  go to $0$ as $n\rightarrow
\infty$ by (2.2) and the condition $\det(\partial^2
u_n),~\det(\partial^2 u)\geq\epsilon_0>0$. The last integral also
goes to $0$ by choosing $\lambda$ small enough. Thus the lemma is
proved.
\end {proof}

Let $u_P$ be  defined  as in (1.1). Then

\begin {lem} For any $u \in {C(\bar{P})}$, $\log^+(\det{(\partial^2 u)})$ is integrable and
 \beq \int_{P}\log^+(\det(\partial^2 u))\leq \int_{P}(u_P)^{ij}_{ij}
udx+\int_{\partial P} u d\sigma+C, \eeq where $C=C(u_P)$.
\end {lem}

\begin {proof} Without loss of generality, we may assume that $u$
satisfies $$\det(\partial^2u), ~\det(D^2u_h)\geq 1,$$ otherwise,
we replace $u$ by $u+\frac{|x|^2}{2}$, and $u_h$ by
$u_h+\frac{|x|^2}{2}$.

For any $\delta > 0$,  we let $P_{\delta}$ be the interior polygon
with faces parallel to those of $P$ separated by distance $\delta$.
 By Lemma 2.2 we have
 \beq \int_{P_{\delta}}\log(\det(\partial^2u))dx= \lim_{h\rightarrow
0}\int_{P_{\delta}}\log(\det(\partial^2u_h))dx.\eeq
 Hence it suffices to prove that
$\int_{P_{\delta}}\log(\det(\partial^2u))dx$ is uniformly bounded
for $\delta$.

By the convexity of $-\log(\det(\cdot))$ for  positive definite
matrices, we have
 $$\log(\det(\partial^2u_h))\leq \log(\det(D^2u_P))+u_P^{ij}(u_h)_{ij}-n.$$
 We need to  estimate $\int_{P_{\delta}}u_P^{ij}(u_h)_{ij}dx$. Using integration  by parts, we obtain
 \beqs &&\int_{P_{\delta}}u_P^{ij}(u_h)_{ij}dx\\
  && = \int_{\partial
P_{\delta}} u_P^{ij}(u_h)_{i}n_{j}d\sigma_0 - \int_{\partial
P_{\delta}}(u_P)^{ij}_{j}u_hn_{i}d\sigma_0+
\int_{P_{\delta}}(u_P)^{ij}_{ij} u_hdx, \eeqs
  where $n$ is the outer normal vector.

Let $x$ be a point of $\partial P_{\delta}$ and $\xi=\xi(x)$ be the
vector $(\xi^i)=(u_P^{ij}n_j)$. Then
$u_P^{ij}(u_h)_{i}n_{j}=\nabla_{\xi} u_h$. Let $y=y(x)$,
$\bar{y}=\bar{y}(x)$ and $z=z(x)$ be the intersection points of the
ray $\{x+t\xi:t>0\}$  with  the boundary $\partial
P_{\frac{\delta}{2}}$, $\partial P_{\frac{3\delta}{2}}$ and
$\partial P$, respectively (see Figure 1 for dimension 2).

\begin{figure}
  \includegraphics[width=10cm]{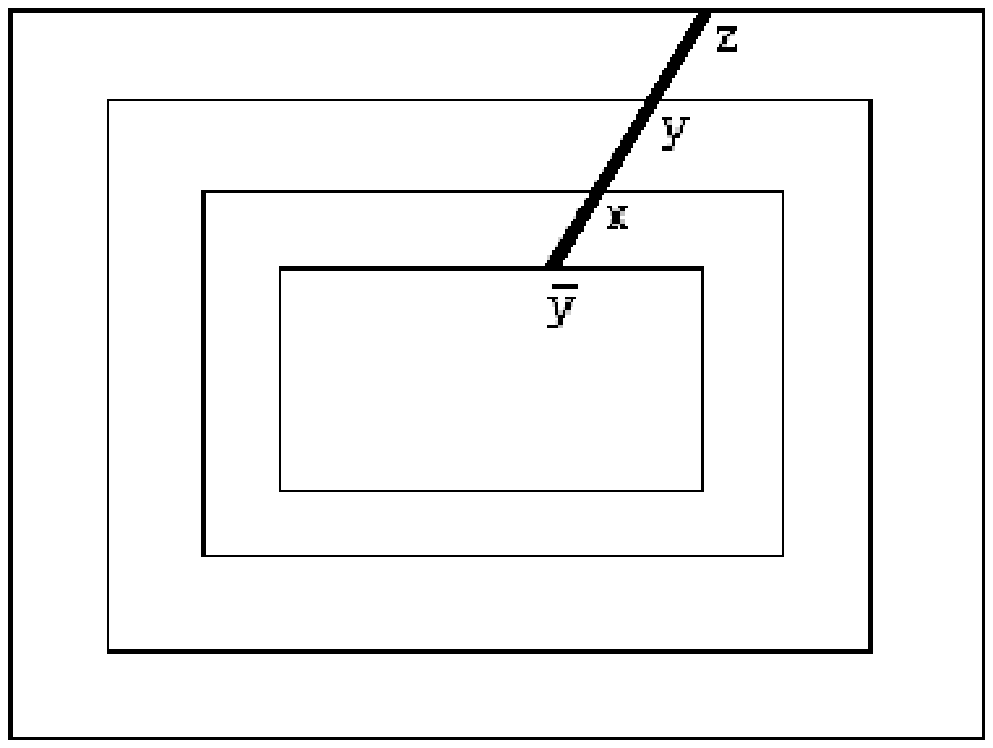}
    \caption{}
\end{figure}
It was verified in [D2] that there exist $c_1,c_2>0$ such that for
any $x\in \partial P_{\delta}$,
 \beq |\xi|\leq
c_1|y-x|=c_1|\bar{y}-x|;~ |\bar{y}-x|=|y-x|\leq c_2\delta. \eeq
 Thus by
the convexity of $u_h$, we have
 \beq |\nabla_{\xi} u_h|\\ \leq
c_1\max\{|u_h(y)-u_h(x)|, |u_h(\bar{y})-u_h(x)|\} \eeq
 and
  \beqn
 &&\int_{\partial P_{\delta}} u_P^{ij}(u_h)_{i}n_{j}d\sigma_0\notag\\
  && \leq c_1
\int_{\partial P_{\delta}}\max\{|u_h(y)-u_h(x)|,
|u_h(\bar{y})-u_h(x)|\}d\sigma_0. \eeqn
 Let $h\rightarrow 0$. By (2.5), we get
 \beqs &&\int_{P_{\delta}}\log(\det(\partial^2u))dx\leq c_1
\int_{\partial
P_{\delta}}\max\{|u(y)-u(x)|, |u(\bar{y})-u(x)|\}d\sigma_0 \\
&&-\int_{\partial P_{\delta}}(u_P)^{ij}_{j}un_{i}d\sigma_0+
\int_{P_{\delta}}(u_P)^{ij}_{ij} udx\\
&&\leq c_1 \int_{\partial
P_{\delta}}(|u(y)-u(x)|+|u(\bar{y})-u(x)|)d\sigma_0 \\
&&-\int_{\partial P_{\delta}}(u_P)^{ij}_{j}un_{i}d\sigma_0+
\int_{P_{\delta}}(u_P)^{ij}_{ij} udx. \eeqs
  Note that the first integral at the last inequality above goes to
$0$ as $\delta\rightarrow 0$ since $u$ is continuous in $\bar{P}$,
and  the second integral  converges to $\int_{\partial P} ud\sigma$
as $\delta$ goes to $0$ since
$(u_P)^{ij}_{j}n_{i}d\sigma_{\delta_0}$ converges to $d\sigma$ on
the boundary $\partial P$ as $\delta\rightarrow 0$ ([D2]). The third
integral clearly converges to $\int_{P}(u_P)^{ij}_{ij} udx$.
Therefore we get (2.3).
\end {proof}

Let us recall an adapted co-ordinates  system on $P$. That is, for
any point $X$ on $\partial P$, there is an affine co-ordinates
$\{x_1,..., x_n\}$, such that there is a neighborhood of $X$ defined
by $p$ inequalities
$$x_1,..., x_p\geq 0,$$
 with  $X$  the origin in this system. Then for  a $v\in \mathcal {C}$, under an adapted
co-ordinates, $v$ has the form
$$v=\sum_{i=1}^{p}x_i\log x_i+w,$$
where $w$ is a smooth function. Furthermore, by choosing a suitable
adapted coordinates system one may assume that at $X=0$,
$$\frac{\partial^2 w}{\partial x_i \partial x_j}=\delta_{ij}, i, j>p,$$
and the converse of $v_{ij}$ has the following property ([D3]),

\begin {lem} \beq(v^{ij})=\text{diag}(f_1x_1, ..., f_px_p,
f_{p+1}, ... ,f_n)+ (\sigma^{ij}), \eeq where $\sigma^{ij}$ is

$$\begin {cases}
&g_{ij}x_ix_j, ~\text{for}~ i, j\leq p,\\
&g_{ij}x_i,  ~\text{for}~ i\leq p, j>p,\\
&h_{ij}, ~\text{for}~ i, j>p,
\end {cases}$$
 and  $f_i$, $g_{ij}$, and $h_{ij}$ are all smooth  functions
 with $f_i(0)=1$ and $h_{ij}(0)=0$ for all $i, j$.
\end {lem}

With this lemma, we are able to deal with functions in $\mathcal
{C}_{\star}$.

\begin {lem} Assume   that $P$ is a $n-$rectangle. Then for any $u \in \mathcal
{C}_{\star}$, $\log^+(\det{(\partial^2 u)})$ is integrable and
 \beq
\int_{P}\log^+(\det(\partial^2 u))\leq \int_{P}(u_P)^{ij}_{ij}
udx+\int_{\partial P} u d\sigma+C, \eeq where $C=C(u_P)$.
\end {lem}

\begin {proof} We still use the notations in Lemma 2.3. As in the
proof of Lemma 2.3, it suffices to prove
$$\int_{\partial
P_{\delta}}(|u(y)-u(x)|+|u(\bar{y})-u(x)|)d\sigma_0\to 0,
~\text{as}~\delta\to 0.$$
 Note that here we could not use the continuity of $u$ up to the boundary $P$.
 We need to verify
  $$\int_{\partial P_{\delta}}|u(y)-u(x)|d\sigma_0\rightarrow 0,
\int_{\partial P_{\delta}}|u(\bar{y})-u(x)|d\sigma_0\rightarrow 0.$$
For simplicity, we just consider the first integral.

 Let  $n-$rectangle $P$ be bounded with $(n-1)$-dimensional faces
 which are defined by
  $$E^i: x_i=\alpha_i, ~ E^{n+i}: x_i=\beta_i\text{~for~}  0<i\leq n,$$
where $\beta_i>\alpha_i$. Then
$$u_P=\sum(x_i-\alpha_i)\log(x_i-\alpha_i)+\sum(-x_i+\beta_i)\log(-x_i+\beta_i).$$
 Let $\{E^i_\delta\}_{i=1}^{2n}$ be the union of corresponding
$(n-1)$-dimensional faces on $\partial P_\delta$.
  We first claim that when $\delta$ is small enough, $z$ must lie on
$E^k$ for any $x$ on $E^k_{\delta}$. This implies that $y$ lies on
$E^k_{\frac{\delta}{2}}$ and $\bar{y}$ lies on
$E^k_{\frac{3\delta}{2}}$. In fact, for any $X$ in $E^k$, we will
make the computation in the adapted coordinates $\{\tilde{x}_1, ...,
\tilde{x}_n\}$ at $X$. Note that since $P$ is a $n-$rectangle, the
transformation of coordinates  just consists of the translations and
reflections. Without loss of generality, we assume that $E^k$
corresponds to $\{\tilde{x}_1=0\}$ in an adapted coordinates. Then
$x=(\delta, \tilde{x}_2, ..., \tilde{x}_n)$ and the outer normal
vector $n=(-1, 0, ...., 0 )$ at $x$.  Since $\tilde{x}_2, ...,
\tilde{x}_p>0$, by Lemma 2.3,  we compute the coordinate of $z$ as
follow
$$\left(0, \tilde{x}_2-\frac{\tilde{x}_2g_{12}\delta}{f_1+g_{11}\delta}, ...,
\tilde{x}_p-\frac{\tilde{x}_pg_{1p}\delta}{f_1+g_{11}\delta},
\tilde{x}_{p+1}-\frac{g_{1,p+1}\delta}{f_1+g_{11}\delta}, ...,
x_{n}-\frac{g_{1,n}\delta}{f_1+g_{11}\delta}\right).$$
 This shows
$\tilde{z}_2, ..., \tilde{z}_p>0$ if $\delta$ is small enough, which
means $z$ lies on $E^k$, so the claim is true in the neighborhood of
$X$. By using the finite covering of adapted coordinates to
$\partial P$,  one can get a uniform small  $\delta$ so that the
claim is true under these adapted coordinates.  Hence the claim is
verified.

For any point $X$ on an open $p$-dimensional face, it lies on the
intersection of $n-p$ faces of dimension $n-1$, denoted by
$E^{k_1}$, ..., $E^{k_{n-p}}$. One chooses a rectangular-shaped
neighborhood $N_{X,\eta}$ in $\bar{P}$ corresponding to a
$n-$rectangle $R_{X,\eta}$
 $$\begin {cases}
&0<\tilde{x}_i<\eta, \text{~for~}  0<i\leq n-p;\\
&-\eta<\tilde{x}_i<\eta, \text{~for~}  n-p+1 \leq i\leq n
\end{cases}$$
in the adapted coordinates of $X$ with $E^{k_j}\bigcap N_{X,\eta}$
corresponding the face $\tilde{x}_j=0$ for any $j\leq n-p$.  In
particular, for a vertex $X$ of $P$, $N_{X,\eta}$ corresponds to
$$0<\tilde{x}_i<\eta, \text{~for~}  0<i\leq n$$
in the adapted co-ordinates. Take the union of all $N_{X,\eta}$
with $X$ on the closed $(n-2)$-dimensional faces and denote it by
$$N_{\eta}=\bigcup_X N_{X,\eta}$$
(see Figure 2, the dark area gives a description of $N_{\eta}$ when
$n=2$ and in higher dimension, it is a neighborhood of all
$(n-2)$-dimensional faces).

\begin{figure}
  \includegraphics[width=7cm]{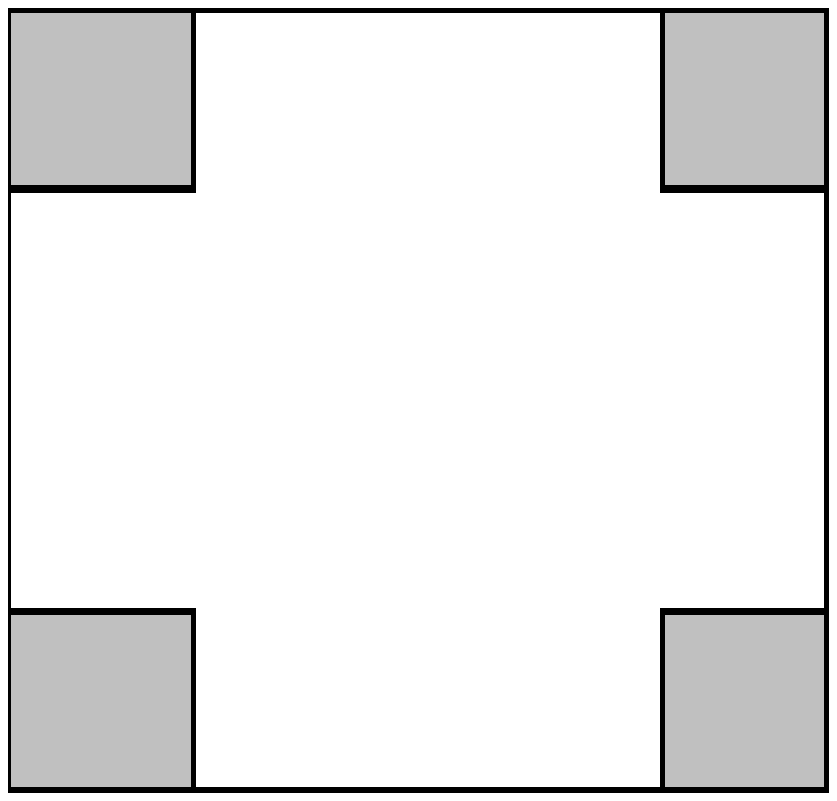}
  \caption{}
\end{figure}

For  a fixed $\eta>0$,  we have
 \beqn &&\int_{\partial P_{\delta}}|u(y)-u(x)|d\sigma_0\notag\\
  &&=\int_{N_{\eta}\cap\partial P_{\delta}}|u(y)-u(x)|d\sigma_0+
\int_{\partial P_{\delta}\setminus N_{\eta}}|u(y)-u(x)|d\sigma_0
\eeqn
  as so as $\delta$ is small. Since we have $|y-x|\leq c_2\delta$ and $|z-x|\leq 2c_2\delta$,
one can suppose $\delta$ to be small enough such that both $y$ and
$z$ lie in
$$\begin {cases}
&N_{2\eta}, ~\text{for}~ x\in{N_{\eta}\cap\partial P_{\delta}};\\
&\bar{P}\setminus N_{\frac{\eta}{2}}, ~\text{for}~ x\in{\partial
P_{\delta}\setminus N_{\eta}}.\end {cases}$$
 Thus  letting
$\delta\rightarrow 0$,  we see that the second integral at the
right hand  of (2.10) goes to zero by the continuity of $u$ in
$\bar{P}\setminus N_{\frac{\eta}{2}}$. It remains to deal with the
first integral. According to the  above computation, for any
$x=(\delta, \tilde{x}_2, ..., \tilde{x}_n)$ in any codimension$-1$
face $E^k_{\delta}$, it must lie in a $N_{X, \eta}$ for some $X
\in {E^k}$. Moreover, in the adapted coordinates, $z$ is a form of
 \beq z(x)=(0, \tilde{x}_2, ..., \tilde{x}_n)+\delta(0, \phi_2(\tilde{x}),
..., \phi_n(\tilde{x})), \eeq
 where $\phi_2(\tilde{x}), ...,
\phi_n(\tilde{x})$ are smooth. This implies that $z$ corresponds
$1-1$ to $x$ and y in $N_{\eta}$ and the measure
$$d\sigma(z)\geq C d\sigma_0(x).$$
Hence by  the convexity of $u$ we obtain
$$\int_{N_{\eta}\cap\partial
P_{\delta}}|u(y)-u(x)|d\sigma_0\leq C \int_{N_{2\eta}\cap\partial
P}u(z)d\sigma.$$
 The late goes to $0$ as $\eta\rightarrow 0$ since
$u$ is integrable on the boundary $\partial P$.
\end {proof}

Let us begin to  prove Proposition 2.1.

\begin {proof}[Proof of Proposition 2.1]
 For any $X$ on $\partial P$, we choose a polytope-shaped
neighborhood $N_X$ in $\bar{P}$ which corresponds to an
$n-$rectangle $R_X$  with an adapted coordinates system
$\{\tilde{x}_1, ..., \tilde{x}_n\}$ at $X$. Since the determinant of
Jacobi matrix associated to the coordinates transformation  is equal
to $1$, we have
$$\int_{N_X}\log(\det(\partial_x^2u))dx=\int_{R_X}\log(\det(\partial_{\tilde{x}}^2u))d\tilde{x}.$$
 Applying Lemma 2.5 to the $n-$rectangle $R_X$ with the corresponding choice of $u_{R_X}$ and using the
 the convexity of $u$,  we obtain
  \beqn \int_{R_X}\log^+(\det(\partial_{\tilde{x}}^2u))d\tilde{x} &&\leq
\int_{R_X}(u_{R_X})^{ij}_{ij} u d\tilde{x}+\int_{\partial R_X} u d\sigma(\tilde{x})+C\notag\\
&&\leq C_1\int_{P}u dx+C_2\int_{\partial P} u d\sigma(x)+C_3, \eeqn
where $C_i=C_i(N_X)$, $i=1,2,3$.

\begin{figure}
  \includegraphics[width=7cm]{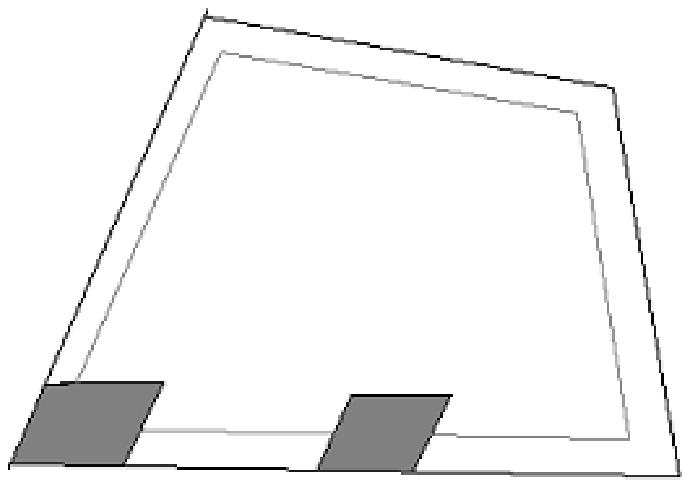}
  \caption{}
\end{figure}

Note that $P\setminus P_\delta$ can be covered by finite
neighborhood $\{N_{X_1}, ..., N_{X_m} \}$ as so as $\delta$ is small
enough (see Figure 3). On the other hand, since $u\in
{C(\bar{P_\delta})}$,  by Lemma 2.3, we have
 \beqn \int_{P_\delta}\log^+(\det(\partial^2 u)) &&\leq
\int_{P_\delta}(u_{P_\delta})^{ij}_{ij} u dx+\int_{\partial
P_\delta} u
d\sigma_\delta+C\notag\\
&&\leq C_1(\delta)\int_{P}u dx+C_2(\delta)\int_{\partial P} u
d\sigma+C_3(\delta), \eeqn
 where $u_{P_\delta}$ is given by
$$\sum (-\langle l_i, x\rangle +
\lambda_i-|l_i|\delta)\log (-\langle l_i, x\rangle +
\lambda_i-|l_i|\delta).$$
 Hence combining (2.12) and (2.13), we get
 \beqs \int_{P}\log^+(\det(\partial^2 u))&&\leq
\sum\int_{N_{X_i}}\log^+(\det(\partial^2u))dx+\int_{P_\delta}\log^+(\det(\partial^2
u))dx\\
&&\leq C_1'\int_{P}u dx+C_2'\int_{\partial P} u d\sigma+C_3', \eeqs
where $C_i'$ are independent of $u$.
\end {proof}

 \begin{rem} By Lemma 2.3 and 2.5,   we also have the following estimate,
 $$\int_{P}\log(\det(\partial^2 u))\leq
r\int_{P}(u_P)^{ij}_{ij} udx+r\int_{\partial P} u d\sigma+C(u_P,
r),~\forall r>0,$$
  since we can replace $u_P$ by $r^{-1}u_P$. Hence  according to the proof of Proposition 2.1,
one   sees that for any small $c$ there exist two uniform $C=C(c)$
and $C'=C'(c)$ such that
 \beq \int_{P}\log(\det(\partial^2 u))\leq c\int_{\partial P} u
d\sigma+C\int_{P}u dx+C'.\eeq
\end{rem}

 By Proposition 2.1, we  get the following approximation result.

\begin {prop}Let $u\in {\mathcal {C}_{\star}^K(P)}$. Suppose that \beq \int_P
\log(\det(\partial^2u))dx>-\infty. \eeq
  Then there exists a sequence
$u_n$ in $\mathcal {C}$ which locally uniformly converges to $u$ and
\beq \int_P  \log(\det(\partial^2u))dx=\lim_{n\rightarrow
\infty}\int_P \log(\det(\partial^2u_n))dx. \eeq
 \end {prop}

\begin {proof} We  prove the proposition by the following three steps.
Note by (2.15) and Proposition 2.1, we see that
$\log(\det(\partial^2u))$ is integrable.

 Step 1. For any
$u\in{\mathcal {C}_{\star}^K(P)}$, there exists a sequence $\{u_n\}$
in $C(\bar{P})$, such that (2.16) holds. In fact we choose
$$u_n(x)=u(r_nx),$$
 with  $r_n\rightarrow1$.  By the integrability of
 $\log(\det(\partial^2u))$, one sees that (2.16) is true.

Step 2. For any $u\in{C(\bar{P})}$, there exists a sequence
$\{u_n\}$ in $\mathcal {C}_{\infty}$, such that (2.16) holds. To
prove this,  we  extend $u$ to a  polytope neighborhood
$P_{-\delta}$ of $P$ with each $n-1$-dimensional face parallel to
one  of  $P$. Then the mollification function $u_h$ is well-define
on $P$ for sufficiently  small $h$.  By  the integrability of
 $\log(\det(\partial^2u))$, it is easy to see
$$\int_P  \log(\det(\partial^2u))dx=\lim_{\epsilon\rightarrow 0}\int_P
\log(\det(\partial^2[u+\epsilon|x|^2]))dx.$$
 On the other hand, since  $P$ can be
regarded as a subset of $P_{-\delta}$, by Lemma 2.2, we have
$$\int_P
\log(\det(\partial^2[u+\epsilon|x|^2]))dx=\lim_{h\rightarrow
0}\int_P \log(\det(\partial^2[u_h+\epsilon|x|^2]))dx.$$
 Thus  the above two relations show that  sequence $u_h+\epsilon|x|^2$ satisfy (2.16).

Step 3. For any $u\in{\mathcal {C}_{\infty}}$, there exists a
sequence $\{u_n\}$ in $\mathcal {C}$, such that (2.16) holds. This
was in fact proved in [D2] where  a sequence was constructed as
follow:  Let $\eta_\delta$ be  a function defined on an interval
 $[-\text{diam}(P),\text{diam}(P)]$, which satisfies
$$\aligned &\eta_{\delta}=x\log x,~\text{if}~ x<\delta,\\
 &\eta''_{\delta}\geq0,~and \\
 & 0\geq \eta_{\delta} \geq 2\delta\log \delta.\endaligned$$
  Let
 $$U_{\delta}(x)=\sum_i\eta_{\delta}(\lambda_i-\langle x, l_i\rangle).$$
  Then the sequence  $u(r_n x)+U_{\delta_n}(x)$ satisfies (2.16).

Combining  the above three steps, we will get  a sequence $u_n$ in
$\mathcal {C}$ such that (2.16) is satisfied.

\end {proof}

From the proof of the above proposition  we see that the sequence
$u_n$ in $\mathcal {C}$ constructed for  $u\in {\mathcal
{C}_{\star}^K(P)}$ satisfying (2.15) also have the property,
$$L(u)=\lim_{n\rightarrow \infty}L(u_n),$$
where
$$L(u)=\int_{\partial P}u
d\sigma-\int_P (\bar{R}+\theta_X)udx$$
 is the linear part of $\mathcal {F}(u)$. Thus  we get

\begin {cor} There exists a $K_0>0$ such that for any $K\ge K_0$ it
holds
  $$\inf_{\mathcal {C}_{\star}^K(P)}\mathcal {F}(\cdot
)=\inf_{\mathcal {C}_{\star}}\mathcal {F}(\cdot )=\inf_{\mathcal
{C}}\mathcal {F}(\cdot).$$
\end {cor}

\section {  Semi-continuity of $\mathcal {F}(u)$}
  \vskip3mm

In this section, we discuss the lower semi-continuity of the
functional $\mathcal {F}(u)$. First we have

\begin {lem}Suppose that ${u_n}\in{\mathcal {C}_{\star}^K(P)}$
converge locally uniformly to $u\in{\mathcal {C}_{\star}^K(P)}$ for
some $K>0$, and \beq \int_P\log(\det(\partial^2u_n))> -C_0 \eeq for
some $C_0>0$. Then for any $h>0$,
$$\limsup_{n\rightarrow\infty}\int_{P_h}\log(\det(\partial^2u_n))\leq\int_{P_h}\log(\det(\partial^2u)).$$
\end {lem}

\begin {proof}By (3.1) and Proposition 2.7, it suffices to prove it
for $u_n \in C^2(P_h)$.  Recall that a convex function on $P$
induces a Monge-Ampere measure $\mu[u]$ through its normal mapping
and this is a Radon measure and can be decomposed into a regular
part and a singular part as follows,
$$\mu[u]=\mu_r[u]+\mu_s[u].$$
Denote by $S$ the supporting set of $\mu_s[u]$, whose Lebesgue
measure is zero. Since ${u_n}$ converges uniformly to $u$, by the
upper semi-continuity of $\mu[u]$,  then for any closed subset
$F\subset P_h\setminus S$,
 \beq \limsup_{n\rightarrow \infty}\int_F
\det(D^2u_n)dx\leq \int_F \det(\partial^2u)dx. \eeq

For given $\epsilon,\epsilon'>0$,  we let
$$\Omega_k=\{x\in P_h\setminus S|~(k-1)\epsilon\leq \log(\det(\partial^2u))<k\epsilon\},
 k=0, \pm1, \pm2,..., $$
and $\omega_k\subset \Omega_k$ be a closed set such that
$|\Omega_k\backslash \omega_k|<\frac{\epsilon'}{2^{|k|}}$. In
particular, we let
$$\Omega_{-\infty}=\{x\in P_h\setminus S|~\log(\det(\partial^2u))=-\infty\}
=\{x\in P_h\setminus S|~\det(\partial^2u)=0\},$$
 and
$\omega_{-\infty}\subset \Omega_{-\infty}$ be a closed set such
that $|\Omega_{-\infty}\backslash \omega_{-\infty}|<\epsilon'$.
Note that $\Omega_{-\infty}$ is a Lebesgue zero set when
$\int_{P_h}\log(\det(\partial^2u))> -\infty$.

First we consider the case that $\Omega_{-\infty}$ is a Lebesgue
zero set. Then for each $\omega_k$, by convexity of $\log$ and
(3.2), we have
 \beqs &&\limsup_{n\rightarrow
\infty}\frac{1}{|\omega_k|}\int_{\omega_k} \log(\det(D^2u_n))dx\\
&&\leq \limsup_{n\rightarrow \infty}\log\left(\frac{\int_{\omega_k}
\det(D^2u_n)dx}{|\omega_k|}\right)\\
&&\leq \limsup_{n\rightarrow
\infty}\log\left(\frac{\int_{\omega_k}
\det(\partial^2u)dx}{|\omega_k|}\right)\\
&&\leq \limsup_{n\rightarrow
\infty}\log\left(\frac{\int_{\omega_k}e^{k\epsilon}
dx}{|\omega_k|}\right)\\
&&\leq k\epsilon. \eeqs
 It follows
  \beqs &&\limsup_{n\rightarrow
\infty}\int_{\omega_k} \log(\det(D^2u_n))dx\\
 && \leq k\epsilon|\omega_k|\\
&&\leq (k-1)\epsilon|\Omega_k|+\epsilon|\Omega_k|+\frac{|k|\epsilon\epsilon'}{2^{|k|}}\\
&&\leq \int_{\Omega_k}
\log(\det(\partial^2u))dx+\epsilon|\Omega_k|+\frac{|k|\epsilon\epsilon'}{2^{|k|}}.
\eeqs
 Hence,
 \beqs &&\limsup_{n\rightarrow\infty}\int_{\bigcup\omega_k}
\log(\det(\partial^2u_n))dx\\ &&\leq \int_{\bigcup\Omega_k}
\log(\det(\partial^2u))dx+ \epsilon
|\bigcup\Omega_k|+C\epsilon\epsilon'\\
&&\leq \int_{\bigcup\Omega_k} \log(\det(\partial^2u))dx+\epsilon
|P_h|+C\epsilon\epsilon'. \eeqs
 Since $\text{osc} u_n$ are uniformly bounded in $P_h$, by (2.2),
 we have
$$\int_{P_h\setminus\bigcup\omega_k} \log(\det(\partial^2u_n))dx\le
C\epsilon'.$$
 Combining the above two inequalities and letting $\epsilon'\rightarrow 0$, we get
$$\limsup_{n\rightarrow \infty}\int_{P_h}\log(\det(D^2u_n))dx\leq
\int_{P_h\setminus S}\log(\det(\partial^2u))dx+C\epsilon.$$
 Letting $\epsilon\rightarrow 0$ again, we obtain
$$\limsup_{n\rightarrow \infty}\int_{P_h}\log(\det(D^2u_n))dx\leq \int_{P_h}\log(\det(\partial^2u))dx.$$

Next we consider the case of $|\Omega_{-\infty}|\neq 0$. In this
case, it must hold
$$\int_{P_h}\log(\det(\partial^2u))=-\infty.$$
By the semi-continuity of $\mu[u]$ one sees
 \beqs
&&\limsup_{n\rightarrow \infty}\frac{\int_{\omega_{-\infty}}
\log(\det(D^2u_n))dx}{|\omega_{-\infty}|} \\ &&\leq
\limsup_{n\rightarrow \infty}\log\left(\frac{\int_{\omega_{-\infty}}
\det(D^2u_n)dx}{|\omega_{-\infty}|}\right)\\
&&\leq \log\left(\frac{\int_{\omega_{-\infty}}
\det(\partial^2u)dx}{|\omega_{-\infty}|}\right)\\
&&= -\infty, \eeqs
  which is contradict to the assumption (3.1).
Thus $|\Omega_{-\infty}|\neq 0$ is impossible. The proof is
finished.
\end {proof}

\begin {lem}Assume that $P$ is a $n-$rectangle. Suppose that
${u_n}\in{\mathcal {C}_{\star}^K(P)}$ converges locally uniformly to
$u\in{\mathcal {C}_{\star}^K(P)}$ for some $K>0$. Then there exists
a subsequence of $\{u_n\}$, still denoted by $\{u_n\}$, and  for any
$\epsilon>0$, there exists  $\delta_0>0$ and a large number $N_0$,
such that for any $0<\delta<\delta_0$ and $n>N_0$, it holds
\beq\int_{P\setminus P_{\delta}}\log(\det(\partial^2u_n))<C\epsilon.
\eeq
 for some uniform constant $C$.
\end {lem}

\begin {proof} Without  loss of generality we may assume that
$$\det(\partial^2u_n)\geq 1.$$
As in Section 2,  we see that there exists a subsequence, still
denoted by $u_n$ converging uniformly on any compact subsets of
codimension$-1$ faces to $\hat{u}$  which satisfies
$$u|_{\partial P}\leq\hat{u}, ~
\int_{\partial P}\hat{u} d\sigma_0=\lim\int_{\partial
P}u_nd\sigma\leq K.$$
 By using  the notations in Section 2,  for
any $\epsilon>0$, we choose $\eta>0$ small enough, such that
$$\int_{N_{2\eta}\cap\partial P}\hat{u}(x)d\sigma<\epsilon.$$

Let $\delta'$ and $\delta$ be two small numbers with
$0<\delta'<\delta$. We need to estimate  the integral for functions
$\log(\det(D^2u_{n}))$   on $P_{\delta'}\setminus P_{\delta}$. Let
$u_{n,h}$ be mollification functions of $u_n$.
 Then
   \beqn
&&\int_{P_{\delta'}\setminus
P_{\delta}}\log(\det(D^2u_{n,h}))\notag \\
 &&\leq\int_{P_{\delta'}\setminus
P_{\delta}}(\log(\det(D^2u_P))-n)+\int_{P_{\delta'}\setminus
P_{\delta}}u_P^{ij}(u_{n,h})_{ij}. \eeqn
 The first integral in the right side of the above inequality  can be made arbitrary
 small as so as $\delta$ is small.  So it suffices to consider the
second integral. Note that using the integration  by parts, we have
 \beqs
&&\int_{P_{\delta'}\setminus
P_{\delta}}u_P^{ij}(u_{n,h})_{ij}\\&&=\int_{P_{\delta'}\setminus
P_{\delta}}(u_P)^{ij}_{ij}u_{n,h}+ \int_{\partial
P_{\delta'}}\nabla_{\xi'} u_{n,h}d\sigma_0-\int_{\partial
P_{\delta}}\nabla_{\xi} u_{n,h}d\sigma_0\\
&&-(\int_{\partial
P_{\delta'}}(u_P)^{ij}_{j}u_{n,h}n_{i}d\sigma_0-\int_{\partial
P_{\delta}}(u_P)^{ij}_{j}u_{n,h}n_{i}d\sigma_0). \eeqs

Let  $x\in{\partial P_{\delta}}$ and $x'\in{\partial
P_{\delta'}}$. Let $y=y(x)$, $\bar{y}=\bar{y}(x)$ and $z=z(x)$ be
the closest points to $x$ on the intersection of the ray
$\{x+t\xi:t>0\}$ with the boundary $\partial
P_{\frac{\delta}{2}}$, $\partial P_{\frac{3\delta}{2}}$ and
$\partial P$,  and $y'=y'(x')$, $\bar{y'}=\bar{y'}(x')$ and
$z'=z'(x')$ are the closest points to $x'$ on the intersection of
the ray $\{x'+t\xi':t>0\}$ with the boundary $\partial
P_{\frac{\delta'}{2}}$, $\partial P_{\frac{3\delta'}{2}}$ and
$\partial P$, respectively. By using the argument in Lemma 2.5, we
have
 \beqn &&|\int_{\partial
P_{\delta'}}\nabla_{\xi'} u_{n,h}d\sigma_0|\notag\\
&&\leq C\int_{\partial P_{\delta'}}
\max\{|u_{n,h}(\bar{y'})-u_{n,h}(x')|,
|u_{n,h}(y')-u_{n,h}(x')|\}d\sigma_0, \eeqn
 and
   \beqn && |\int_{\partial
P_{\delta}}\nabla_{\xi} u_{n,h}d\sigma_0|\\\notag
 &&\leq
C\int_{\partial P_{\delta}} \max\{|u_{n,h}(\bar{y})-u_{n,h}(x)|,
|u_{n,h}(y)-u_{n,h}(x)|\}d\sigma_0, \eeqn
 Letting $h\rightarrow 0$ in both (3.5) and (3.6), then we get
 \beqn &&\int_{P_{\delta'}\setminus
P_{\delta}}\log(\det(\partial^2u_n))\notag\\
&&\leq\int_{P_{\delta'}\setminus
P_{\delta}}(u_P)^{ij}_{ij}u_n-(\int_{\partial
P_{\delta'}}(u_P)^{ij}_{j}u_nn_{i}d\sigma_0-\int_{\partial
P_{\delta}}(u_P)^{ij}_{j}u_nn_{i}d\sigma_0)\notag\\
 &&+C\int_{\partial
P_{\delta'}} (|u_n(\bar{y'})-u_n(x')|+
|u_n(y')-u_n(x')|)d\sigma_0\notag \\
&& + C \int_{\partial P_{\delta}} (|u_n(\bar{y})-u_n(x)|+
|u_n(y)-u_n(x)|)d\sigma_0. \eeqn

As in Lemma 2.5, we shall use the adapted coordinates on a union of
several $n$-rectangles $N_{\eta}$. By the convergence of $u_n$ to
$u$, we see that there exists $N_0$ such that for $n>N_0$ and $x\in
\bar{P}\setminus N_{\frac{\eta}{2}}$,
$$|u_n(x)-u(x)|<\epsilon.$$
Also by the  continuity of $u$ in $\bar{P}\setminus
N_{\frac{\eta}{2}}$, for any $x$ in $\bar{P}\setminus N_{\eta}$, we
have
$$|u(y)-u(x)|, ~|u(\bar{y})-u(x)|, ~|u(y')-u(x')|, ~|u(\bar{y'})-u(x')|<\epsilon.$$
 Then letting $\delta$ be small enough, we get the following
estimates
 \beq \int_{(P_{\delta'}\setminus
P_{\delta})\bigcap(\bar{P}\setminus N_{\eta})}(u_P)^{ij}_{ij}u_n\leq
C\int_{P_{\delta'}\setminus P_{\delta}}u+C\epsilon\leq C'\epsilon,
\eeq
 \beq \int_{\partial P_{\delta'}\bigcap(\bar{P}\setminus
N_{\eta})} (|u_n(\bar{y'})-u_n(x')|+ |u_n(y')-u_n(x')|)d\sigma_0\leq
C\epsilon, \eeq
 and
   \beq \int_{\partial
P_{\delta}\bigcap(\bar{P}\setminus N_{\eta})}
(|u_n(\bar{y})-u_n(x)|+ |u_n(y)-u_n(x)|)d\sigma_0\leq C\epsilon.
\eeq
 Note that
$(u_P)^{ij}_{j}n_{i}d\sigma_0$ converge to $d\sigma$ on the boundary
$\partial P$ as $\delta$ goes to $0$. Thus as so as $\delta$ is
small enough we obtain
  \beqs &&|\int_{\partial
P_{\delta'}\bigcap(\bar{P}\setminus
N_{\eta})}(u_P)^{ij}_{j}u_nn_{i}d\sigma_0-\int_{\partial
P_{\delta}\bigcap(\bar{P}\setminus
N_{\eta})}(u_P)^{ij}_{j}u_nn_{i}d\sigma_0|\\
&&\leq |\int_{\partial P_{\delta'}\bigcap(\bar{P}\setminus
N_{\eta})}(u_P)^{ij}_{j}un_{i}d\sigma_0-\int_{\partial
P_{\delta}\bigcap(\bar{P}\setminus
N_{\eta})}(u_P)^{ij}_{j}un_{i}d\sigma_0|+C\epsilon\\
&&\leq C\epsilon. \eeqs

For the integrals on $N_{\eta}$, we use the convexities of $u_n$ and
estimate
  \beq \int_{(P_{\delta'}\setminus P_{\delta})\bigcap
N_{\eta}}(u_P)^{ij}_{ij}u_n\leq C\int_{\partial P\bigcap
N_{2\eta}}u_nd\sigma, \eeq
 \beqn &&\int_{\partial P_{\delta'}\bigcap
N_{\eta}} (|u_n(\bar{y'})-u_n(x')|+
|u_n(y')-u_n(x')|)d\sigma_0\notag \\
&&\leq C\int_{\partial P\bigcap N_{2\eta}}u_nd\sigma, \eeqn
 \beqn
&&\int_{\partial P_{\delta}\bigcap N_{\eta}} (|u_n(\bar{y})-u_n(x)|+
|u_n(y)-u_n(x)|)d\sigma_0\notag\\
&&\leq C\int_{\partial P\bigcap N_{2\eta}}u_nd\sigma, \eeqn
 \beq
\int_{\partial P_{\delta}\bigcap
N_{\eta}}(u_P)^{ij}_{j}u_nn_{i}d\sigma_0\leq\int_{\partial P\bigcap
N_{2\eta}}u_nd\sigma+C\epsilon, \eeq
 and
\beq \int_{\partial P_{\delta}\bigcap
N_{\eta}}(u_P)^{ij}_{j}u_nn_{i}d\sigma_0\leq\int_{\partial P\bigcap
N_{2\eta}}u_nd\sigma+C\epsilon. \eeq
 Note that $u_n$ converges
uniformly on any compact subsets of $\partial P\bigcap N_{2\eta}$.
By choosing $n$ sufficiently large, we have
  \beq \int_{\partial
P\bigcap N_{2\eta}}u_nd\sigma\leq \int_{\partial P\bigcap
N_{2\eta}}\hat{u}d\sigma+C\epsilon\leq C\epsilon \eeq
 Therefore combining
(3.7)-(3.16),  we finally get from (3.4),
$$\int_{P_{\delta'}\setminus
P_{\delta}}\log(\det(\partial^2u_n))\leq C\epsilon.$$
 By letting $\delta'\rightarrow 0$, we obtain (3.3).
\end {proof}

By using the above lemma and the argument in the proof of
Proposition 2.1, we can generalize Lemma 3.2 to the general
polytope $P$ which satisfies the delzant's condition.

\begin {lem} Let  $P$ be a polytope which satisfies the delzant's condition.
 Suppose that
${u_n}\in{\mathcal {C}_{\star}^K(P)}$ converges locally uniformly to
$u\in{\mathcal {C}_{\star}^K(P)}$ for some $K>0$. Then there exists
a subsequence of $\{u_n\}$, still denoted by $\{u_n\}$, and for any
$\epsilon>0$, there exist $\delta_0>0$ and a large number $N_0$ such
that for any $0<\delta<\delta_0$ and $n>N_0$, (3.3)  holds.
\end {lem}

\begin {proof} We also assume that
$$\det(\partial^2u_n)\geq 1.$$
As in the proof Proposition 2.1,  for a   small enough $\delta$, we
choose a covering of  finite polytope-shaped neighborhoods
$\{N_{X_i}\}_{i=1}^m$ of $P\setminus P_{\delta}$ which correspond to
$n-$rectangles $\{R_{X_i}\}_{i=1}^m$   with the adapted coordinates.
Then we observe
  \beqs \int_{P\setminus
P_{\delta}}\log(\det(\partial^2u_n))&&\leq \sum
\int_{N_{X_i}\setminus P_{\delta}}\log(\det(\partial^2u_n))\\
&&\leq \sum \int_{N_{X_i}\setminus
(N_{X_i})_{\delta}}\log(\det(\partial^2u_n))\\, &&\leq \sum
\int_{R_{X_i}\setminus
(R_{X_i})_{\delta'}}\log(\det(\partial^2u_n)) \eeqs
  where
$(N_{X_i})_{\delta}$ and $(R_{X_i})_{\delta'}$ are the interior
polygons corresponding to $N_{X_i}$ and $R_{X_i}$ with faces
parallel to those of $R_{X_i}$ separated by distance $\delta$ and
$\delta'$. Note that $\delta'$ is less than a scalar multiple of
$\delta$ by coordinates transformation.
 Since ${u_n}\in{\mathcal
{C}_{\star}^K(P)}$ implies ${u_n}\in{\mathcal
{C}_{\star}^{K_i}(R_{X_i})}$ for some $K_i>0$ by the convexity of
$u$, the lemma   will follow by applying Lemma 3.2 to each
$R_{X_i}$.

\end {proof}

Now we  prove the main result in this section.

\begin {prop}Suppose that ${u_n}\in{\mathcal {C}_{\star}^K(P)}$
converge locally uniformly to $u\in{\mathcal {C}_{\star}^K(P)}$
for some $K>0$, and $u_n$ satisfies (3.1). Then
$$\int_P\log(\det(\partial^2u))>-\infty,$$
and there exists a subsequence  of  $u_n$ such that
  \beq
\limsup_{n\rightarrow\infty}\int_{P}\log(\det(\partial^2u_n))\leq\int_{P}\log(\det(\partial^2u))
\eeq
\end {prop}

\begin {proof} On the contrary, we  assume
 \beq \int_P\log(\det(\partial^2u))= -\infty.
\eeq
 Then by Proposition 2.1,  we see that for any $M>0$, there
 exists
$\delta_0>0$ such that for any $\delta<\delta_0$,
 \beq
\int_{P_{\delta}}\log(\det(\partial^2u))<-M. \eeq
 On the other
hand, by Lemma 3.3, there exist $\delta_1>0$ and $N_1>0$ such that
for any $0<\delta<\delta_1$ and $n>N_1$
 \beq \int_{P\setminus
P_{\delta}}\log(\det(\partial^2u_n))<C\epsilon. \eeq
 Then for a fixed $\delta<\min(\delta_0, \delta_1)$, by Lemma 3.1,  we see that
 there exists
$N_2>0$ such that for $n>N_2$,
 \beq
\int_{P_{\delta}}\log(\det(\partial^2u_n))\leq
\int_{P_{\delta}}\log(\det(\partial^2u))+\epsilon. \eeq
 Thus by choosing $n>N_0=\max(N_1, N_2)$ and using  (3.19)-(3.21), we
get
 \beqs
\int_P\log(\det(\partial^2u_n))&&\leq\int_{P_{\delta}}\log(\det(\partial^2u_n))+C\epsilon\\
&&\leq\int_{P_{\delta}}\log(\det(\partial^2u))+(C+1)\epsilon\\
&&\leq -M+(C+1)\epsilon. \eeqs
 Since $M$ is arbitrary, we derive
   $$\limsup_{n\rightarrow\infty}\int_{P}\log(\det(\partial^2u_n))=-\infty.$$
This is a contradiction.  Hence we prove
  \beq -\infty<
\int_P\log\det(\partial^2u)< \infty. \eeq

By (3.22) we see that there exist $\delta_0>0$ and $N_3>0$ such
that for any  $0<\delta<\delta_0$,
 \beq |\int_{P\setminus
P_{\delta}}\log(\det(\partial^2u))|<\epsilon. \eeq
 Combining  (3.20) and (3.21),  we  get
\beqs
\int_P\log(\det(\partial^2u_n))&&\leq\int_{P_{\delta}}\log(\det(\partial^2u_n))+C\epsilon\\
&&\leq\int_{P_{\delta}}\log(\det(\partial^2u))+(C+1)\epsilon\\
&&\leq\int_P\log(\det(\partial^2u))+(C+2)\epsilon. \eeqs
 Letting
$\epsilon\to 0$, we will obtain (3.16)
\end {proof}

\begin{rem} Proposition 3.4 implies that the functional $\mathcal
{F}(u)$ is lower semi-continuous in $\mathcal C_\star$ since $L(u)$
 is lower semi-continue in $\mathcal C_\star$.
 \end {rem}

\section { Proof of  Theorem 0.2}
  \vskip3mm

In this section, we prove the existence of minimizing weak solution
for the extremal metrics.

\begin {prop} Suppose that  the modified $K$-energy $\mu(\phi)$ is proper in
$\mathcal {M}_{G_0}$ associated  to toric actions group $T$ on a
toric manifold $M$. Then for any minimizing sequence $u_n$ in
$\tilde{\mathcal {C}}$
 there exists a subsequence of $u_n$, still denoted by $u_n$, which  locally uniformly converge to
a convex function $u_{\infty}$ such that
$$\mathcal {F}(u_\infty)=\inf_{\tilde{\mathcal {C}}}\mathcal {F}(\cdot).$$
 Moreover, $u_n|_{\partial P}\to u_{\infty}|_{\partial P}$ as $n\to
 \infty$ almost everywhere.

\end {prop}

\begin {proof}
 By the assumption of properness of $\mu(\phi)$ and Proposition 1.3, one sees that
  there exist constant $C_1$ and $C_2>0$ such that for any normalized
minimizing sequence $u_n$ in $\tilde{\mathcal {C}}$,
 \beq
\int_{P}u_n dx <C_1, ~ -\int_P\log(\det(\partial^2u_n))+
\int_{\partial P}u_n d\sigma<C_2. \eeq
 On the other hand, by (2.14) in Remark 2.6, we know that for a positive $c<1$
 there exist $C=C(c)$ and $C'=C'(c)$ such that
  \beq
-\int_P\log(\det(\partial^2u_n))\geq -c\int_{\partial P}u_n
d\sigma-C\int_{\partial P}u_n dx-C'. \eeq
 Combining these two inequalities, we get
  \beq \int_{\partial P}u_nd\sigma\le \max_{i=1,...,d}\{\frac{1}{|l_i|}\}\int_{\partial P}u_nd\sigma_0<K_0 \eeq
 for  some $K_0>0$. This implies  $u_n\in{\mathcal {C}_{\star}^K(P)}$ for any $K\ge K_0$.
 Thus by the convexities of $u_n$ there exists  a subsequence of $\{u_n\}$, still denoted by $\{u_n\}$
which locally uniformly converges  to a normalized function
$u_{\infty}$ in $\mathcal {C}_{\star}(K_0)$ satisfying
 $$\int_{\partial P}u_{\infty}d\sigma\leq\liminf_{n\rightarrow\infty}\int_{\partial P}u_nd\sigma<K.$$
  On the other hand, by (4.3) together with (4.1) and (4.2), we also have
$$\int_P\log(\det(\partial^2u_n))> -C_0$$
 for some $C_0$. Then by Remark 3.5, we get
$$\mathcal {F}(u_{\infty})\leq \liminf_{n\rightarrow\infty}\mathcal {F}(u_n)
=\inf_{\mathcal {\tilde C}} \mathcal {F}(u).$$
 Hence by Corollary 2.8, we prove
$$\mathcal {F}(u_{\infty})=\inf_{\mathcal {C}_{\star}^K(P)}\mathcal {F}(\cdot)
=\inf_{\tilde{\mathcal {C}}}\mathcal {F}(\cdot).$$

It remains to prove
$$u_{\infty}|_{\partial P}=\hat{u},$$
where $\hat{u}=\lim_{n\to\infty} u_n|_{\partial P}$. In fact,
$$-\int_{P}\log(\det(\partial^2u_\infty)) dx+\int_{\partial
P}\hat{u}d\sigma-\int_{P}u_{\infty}dx$$
  is  the minimum of the functional $\mathcal {F}(\cdot)$ in $\mathcal {C}_{\star}^K(P)$. This implies
$$\int_{\partial P}\hat{u}d\sigma=\int_{\partial P}u_{\infty}d\sigma.$$
Thus by the fact $u_{\infty}|_{\partial P}\leq\hat{u}$, we conclude
that $u_n|_{\partial P}\to u_{\infty}|_{\partial P}$ as $n\to
 \infty$ almost everywhere.

\end {proof}

In  the proof of (4.3), we used the estimate (4.2) established in
Section 2. In the following we give a direct proof to (4.3).
 Let $\{u_i\}$ be  a minimizing  sequence of $\mathcal {F}(u)$ in
$\tilde{\mathcal {C}}$. By (1.5), we may assume that for each $i$ it
holds
$$\mathcal {F}(u_i)=\inf_\lambda \mathcal {F}(\lambda  u_i).$$
 Note  $\lambda  u_i\in\mathcal {C}_\infty$ for any $\lambda>0$.
  We claim
   \beq \int_{\partial P} u_i d\sigma=t_i\le K.\eeq
Let $\overline u_i=\frac{u_i}{t_i}$ so that
 $$ \int_{\partial P} \overline u_i d\sigma=1,~\text{for each}~, i=1,2,...$$
   Then $t_i$ is a critical point of function $\mathcal {F}(\lambda  u_i)$ for
   $\lambda$. So we have
$$\frac{n Vol(P)}{t_i}=L(\overline u_i).$$
 Thus the claim  will be true if
 $$L(\overline u_i)\ge \delta>0,~\text{for sufficiently large }~i.$$
   Suppose that the above is not true. Then there exists a subsequence $\overline u_{i_k}$ such that
\beq \lim_{i_k}L(\overline u_{i_k})=0.\eeq
  It follows
$$\lim_{i_k} \int_P(\bar{R}+\theta_X)\overline {u}_{i_k} dx=\lim_{i_k} \int_{\partial P} \overline{u}_{i_k}d\sigma
 = 1.$$
Thus
 $$\int_P \overline {u}_{i_k} dx \ge \delta'>0, ~\text{for sufficiently large }~i.$$
  On the other hand, by
the properness of $ \mathcal {F}(u)$, we have
 $$C\ge \mathcal {F}(u_i)\ge p( t_i\int_P \overline u_i dx), ~\text{for each}~, i=1,2,...$$
Hence  we get $t_{i_k}\le C'$ and
 $$L(\overline u_{i_k})=\frac{n}{t_{i_k}}\ge \frac{n}{C'}.$$
  The late is contradict to the assumption (4.5). Therefore the claim is
  true. (4.4) implies (4.3).


 Theorem 0.2   follows from Proposition 4.1.   We also have the
 following corollary.

\begin {cor} Let $M$ be a toric K\"ahler
manifold associated to a convex polytope $P$ described by (0.3).
Suppose
 \beq \bar{R}+\theta_X < \frac{n+1}{\lambda_i}, ~\forall
~i=1,..., d. \eeq
 Then there exists  a minimizing  weak solution of equation (0.5)
  in $M$ in the sense of convex functions.
\end {cor}

\begin{proof} It was proved in [ZZ] that under the condition  (4.6)
we have
$$\mathcal {F}(u)\ge c\int_P u dx-C,~\forall u\in \mathcal {\tilde {\mathcal
C}}$$
 for some $c, C>0.$ In particular, the modified $K$-energy $\mu(\phi)$  is
proper  in $\mathcal {M}_{G_0}$ associated  to toric actions group
$T$. Thus by Proposition 4.1, the corollary is true.
\end{proof}

Next we discuss some  properties  about  the minimizing weak
solution $u_\infty$.

\begin {prop} Suppose that
$$\theta_X +\overline R\ge 0.$$
 Then for any minimizer $u_\infty$ of $\mathcal {F} (\cdot)$, the
Monge-Ampere measure $\mu[u_\infty]$ has no singular part.
\end {prop}

\begin {proof} We use an argument from [TW2] to prove the proposition.
 Suppose $\mu[u_\infty]$ has non-vanishing singular part
$\mu_s[u_\infty]$. Then  for any $M>0$, there must exist a ball
$B_r\subset P$ such that
 \beq \mu_s[u_\infty](B_r)\geq
M(\mu_r[u_\infty](B_r)+|B_r|). \eeq
 We consider  the following
Dirichlet problem for Monge-Amp\`ere operator,
$$ \begin{cases} &\mu[v]=M\mu_r[u_\infty]+ M ~\text{in} ~ B_r,\\
&v=u_\infty ~\text{on} ~ \partial B_r.\end{cases}$$
 By the Alexander theorem, the above equation has  a unique
convex solution $v$. Note
 \beq
\det(\partial^2v)=M\det(\partial^2 u_\infty)+ M,~\text{in}~ B_r.
\eeq
 By comparison principle, $u_\infty\leq v$ in $B_r$, and  the set
$E=\{v>u_\infty\}$ is not empty. Define another convex  function
$\tilde u$ by
 $$\begin {cases} &\tilde{u}= u ~\text{in}~ P\backslash E,\\
&\tilde{u}=v ~\text{in}~ E.\end {cases}$$
 We claim $\mathcal
{F}(\tilde{u})<\mathcal {F}(u_\infty)$, so we get a contradiction to
the assumption that $u$ is a minimizer.  In fact, by choosing $M$
sufficiently  large and using (4.8),  we have
 \beqs &&\mathcal {F}(\tilde{u})-\mathcal {F}(u_\infty)\\
&&=-\int_E \log(\det (\partial^2v))dx+\int_E \log(\det(\partial^2u_\infty))dx
-\int_E (\theta_X+\overline R)(v-u_\infty)dx\\
&&\leq -(\log M)|E|-\int_E
\log(\det(\partial^2u_\infty)+1)dx+\int_E\log(\det(\partial^2u_\infty))dx\\
&&<0. \eeqs
 The proposition is proved.

\end {proof}

Recall $u$ a piecewise linear (PL) function   on $P$ if $u$ is a
form of
$$u= \text{max}\{u^1,...,u^r\},$$
 where $u^\lambda=\sum a_i^\lambda x_i + c^\lambda, ~\lambda=1,...,r,$
for some vectors $(a^\lambda_i)\in \mathbb {R}^n$ and some numbers
$c^\lambda\in \mathbb {R}$. In particular, when $r=2$ and $u^2=0$,
we call $u$ simple PL function with crease $\{u^1=0\}$. The
following Lemma  was proved in [ZZ].

\begin {lem} Suppose that for each $i=1,...,d$, it holds
$$\bar{R}+\theta_X < \frac{n+1}{\lambda_i}, ~\text{in}~
P.$$
  Then for any PL-function $u$ on $P$, we have
$$L(u)\geq 0.$$
 Moreover the equality  holds if and only if $u$ is  an affine linear function.
\end {lem}

\begin {prop} Suppose that for each $i=1,...,d$, it holds
$$\bar{R}+\theta_X < \frac{n+1}{\lambda_i}, ~\text{in}~
P.$$
 Let $u$ be  a minimizer of $\mathcal {F}(\cdot)$. Then for any point
$x$ in $P$, there exists only one supporting hyperplane at $x$.
\end {prop}

\begin {proof} Suppose by contrary that there are two distinct
supporting hyperplane at $x_0$ in $P$. By adding an affine linear
function, we  may assume that one of them is zero hyperplane and the
other is
 $$\sum_{i=1}^{n}a_i x_i-x_{n+1}+c=0.$$
  Then we define a simple PL function $f$ with crease through $x_0$
  by
$$f=\max\{\sum_{i=1}^{n}a_i x_i+c, 0\}$$
and  write  $u$ as
$$u=f+(u-f).$$
 It is easy to see that  $u-f$  is a convex function $P$.
    Since $x_0$ lies in $P$, by the
lemma above we have $L(f)> 0$.  On the other hand,
$$\det (\partial^2u)= \det (\partial^2(u-f))$$
almost everywhere. Thus
$$\mathcal {F}(u-f)<\mathcal {F}(u),$$
  which is  contradict to that
$u$ is minimizer. Hence the proposition is true.
\end {proof}

\end{document}